\documentclass[a4paper]{article}
\usepackage{graphicx}
\usepackage{amsmath}
\usepackage{amsfonts}
\usepackage{amsthm}
\usepackage{amssymb}
\usepackage{subfigure}
\usepackage{multirow}
\usepackage{rotating}
\usepackage{fancybox}
\usepackage{boxedminipage}
\usepackage{version}
\usepackage[english,greek]{babel}
\usepackage[iso-8859-7]{inputenc}
\newcommand{\ds}{\displaystyle}
\newcommand{\en}{\selectlanguage{english}}

\newcommand{\NN}{\mathbb{N}}

\newcommand{\ZZ}{\mathbb{Z}}

\newcommand{\FF}{\mathbb{F}}

\newtheorem{thm}{Theorem}

\newtheorem{heur}{Heuristic}
\theoremstyle{definition}
\newtheorem{rmk}{Remark}
\newtheorem{dfn}{Definition}

\newtheorem{cor}{Corollary}
\newtheorem{ex}{Example}
\newtheorem{alg}{Algorithm}

\begin{document} \en

\title{On the generalization of the Costas property to higher dimensions}
\author{Konstantinos Drakakis\footnote{The author holds a Diploma in Electrical and Computer Engineering from NTUA, Athens, Greece, and a Ph.D. in Applied and Computational Mathematics from Princeton University, NJ, USA. He was a scholar of the Lilian Boudouris Foundation.}
\\Electronic and Electrical Engineering\\University College Dublin\\ \& \\ Claude Shannon Institute\footnote{www.shannoninstitute.ie}\\Ireland}
\maketitle

\abstract{We investigate the generalization of the Costas property in 3 or more dimensions, and we seek an appropriate definition; the 2 main complications are a) that the number of ``dots'' this multidimensional structure should have is not obvious, and b) that the notion of the multidimensional permutation needs some clarification. After proposing various alternatives for the generalization of the definition of the Costas property, based on the definitions of the Costas property in 1 or 2 dimensions, we also offer some construction methods, the main one of which is based on the idea of reshaping Costas arrays into higher-dimensional entities.}

\section{Introduction}

Costas arrays originated in engineering as a frequency hopping pattern that optimizes the performance of radars and sonars \cite{C1,C2}; being a singular combinatorial object, however, they have been lately the focus of intensive study by mathematicians as well, and have thus started leading a second independent ``life'' in the mathematical literature \cite{G,GT,G2,GT2,M}. In this work, we increase the level of mathematical abstraction by investigating analogs of Costas arrays in 3 or more dimensions: several challenges lie ahead, as we first need to give a satisfactory definition of the Costas property in higher dimensions, and also provide some algorithms to construct such Costas cubes or ``hypercubes'' in general.

We will begin by stating the definition of Costas arrays \cite{D} in such a way as to exhibit their direct relation to Golomb rulers \cite{R}, which will consequently be construed as the 1-dimensional analog of Costas arrays. Subsequently, we will discuss several slightly different candidate definitions in higher dimensions, and we will attempt to generalize the generation methods available to us. We shall see, in particular, that it is possible to generate sparse Costas hypercubes systematically, using various construction methods, but if we want to know what the densest Costas hypercube can be, things become less clear.    

\section{The definition of Costas property in 1 and 2 dimensions}

We use the usual notation that $[n]:=\{1,2,\ldots,n\},\ n\in\NN$, and that $[n]-1:=\{0,1,\ldots,n-1\},\ n\in\NN^*$. 

\begin{dfn}[Costas property in 1 dimension] \label{c1d} Let $g:[N]\rightarrow [N'],\ N,N'\in\NN$, be a finite strictly increasing sequence, such that $g(1)=1$, $g(N)=N'$; $g$ will satisfy the \emph{Costas property in one dimension} iff $g(i)-g(i+k)=g(j)-g(j+k)\Leftrightarrow i=j,\ i,j,i+k,j+k\in[N]$. Equivalently, if $f:[N']\rightarrow \{0,1\}$ is a sequence such that $f(i)=1\Leftrightarrow \exists j\in[N]: i=g(j)$, and $f(i)=0$ whenever $i\notin [N']$, then the \emph{autocorrelation} $A_f(k)=\sum_{i\in[N']}f(i)f(i+k),\ k\in\ZZ$ satisfies $A_f(k)\leq 1,\ k\neq 0$. The sequence $g$ is known a \emph{Golomb ruler} \cite{R}.
\end{dfn}

\begin{dfn}[Costas property in 2 dimensions] \label{c2d} Let $f:\ZZ^2\rightarrow \{0,1\}$ be a binary sequence in 2 dimensions, such that it is equal to 1 at a finite number of points only; without loss of generality we may assume that $f(i,j)=0$ when $i\notin [M]$ or $j\notin[N]$ for some $M,N\in\NN$. Defining the \emph{autocorrelation} as $A_f(k,l)=\sum_{i,j\in\ZZ}f(i,j)f(i+k,j+l),\ k,l\in\ZZ$, $f$ will have the \emph{Costas property in 2 dimensions} iff $A_f(k,l)\leq 1,\ k,l\neq 0$.
\end{dfn}

\begin{dfn}[Costas squares] \label{cs} Let $f$ have the Costas property in 2 dimensions and let $M,N$, appearing in Definition \ref{c2d} be the least possible; if $M=N=n$, $f$ will be called a \emph{Costas square} of side length $n$. 
\end{dfn}

\begin{dfn}[Costas arrays] \label{ca} Let $f$ be a Costas square, and let it, in addition, have the structure of a permutation matrix with only one element equal to 1 per row and column: in other words, let there be a bijective sequence $g:[n]\rightarrow [n]$ (namely a \emph{permutation} of order $n$) such that $f(i,j)=1\Leftrightarrow i=g(j), i,j\in[n]$. In that case, $f$ will be called a \emph{Costas array}, and $g$ a \emph{Costas permutation}.
\end{dfn}

\begin{rmk}\ \label{car} 

\begin{itemize}
	\item It follows immediately that, if $g$ is a Costas  permutation of order $n$, no two vectors in the collection $\{(i-j, g(i)-g(j)):\ 1\leq i<j\leq n\}$ can be equal, and that no vector in this collection has a coordinate equal to 0. 
	\item It is customary to denote 0s in the array $f$ by blanks and 1s by dots. 
	\item Note that the Costas property, as stated in Definition \ref{c2d}, is (much) more general than usual, as it can be satisfied by (many) more 2-dimensional sequences than Costas arrays themselves, which appear only as quite restricted special cases (see Definition \ref{ca}). But this generality will actually prove to be helpful when we seek a higher-dimensional analog of Costas arrays.
\end{itemize}
\end{rmk}

\section{Costas hypercubes}

It is straightforward to generalize the Costas property in higher dimensions, to the extent that we risk (re)stating the obvious below. What is harder is to find the exact higher-dimensional analog of the special case, the ``Costas array''.

\subsection{The Costas property in higher dimensions}

\begin{dfn}[Costas hyper-rectangles] \label{chr}
Let $m\in \NN$, consider a sequence $f: \ZZ^m\rightarrow \{0,1\}$, and suppose further that $f(i)=0,\ i\notin [N], N\in\NN^m$, where $i=(i_1,\ldots,i_m)$, $N=(N_1,\ldots,N_m)$, $[N]=[N_1]\times\ldots\times[N_m]$, and the vector $N$ has the smallest possible entries (for the given sequence $f$). Let the autocorrelation of $f$ be $A_f(k)=\sum_{i\in\ZZ^m}f(i)f(i+k),\ k\in\ZZ^m$. Then, $f$ will be a \emph{Costas hyper-rectangle} iff $\forall k \in \ZZ^m-\{0\},\ A_f(k)\leq 1$.
\end{dfn}

\begin{dfn}[Costas hypercubes]
If $f$ is a Costas hyper-rectangle so that $N_1=\ldots=N_m=n\in\NN$, it is called a \emph{Costas hypercube}; if $m=2$, it is a Costas square (see Definition \ref{cs}), whereas if $m=3$ it will be called a \emph{Costas cube}. 
\end{dfn}
Now it is time to seek the higher-dimensional analog of a permutation.

\subsection{Vector permutations} \label{vec}

It is quite simple to generalize permutations in higher dimensions, if the number of dimensions is even:

\begin{dfn}[Permutation Costas hypercube]
Let $m=2s,\ s\in\NN$, and let $g:[n]^s\rightarrow [n]^s$ be a bijection, that is a permutation on vectors in general. Let $f: \ZZ^m\rightarrow \{0,1\}$ be a sequence such that $f(i)=1$ iff $(i_{s+1},\ldots,i_{2s})= g(i_1,\ldots,i_s),\ (i_1,\ldots,i_s)\in [n]^s$, and such that it has the Costas property, as defined in Definition \ref{chr}; then, $f$ will be called a \emph{permutation Costas hypercube} in $m$ dimensions with side length $n$. 
\end{dfn}

\begin{rmk} \

\begin{itemize}
	\item Permutation Costas hyper-rectangles are defined by the obvious extension of the above definition. 
	\item The fact that we chose the first $s$ dimensions to form the domain of $g$ and the last $s$ its range does not affect generality: if $f$ is a Costas hypercube, then any $f'$ resulting by a random permutation of the order of the dimensions is also a Costas hypercube; this constitutes a generalization of the invariance of the Costas property in 2 dimensions under transposition.
	\item No 2 vectors in the collection $\{(i-j, g(i)-g(j)):\ i,j\in[n]^s, i\neq j\}$ can be equal; however, they may have coordinates equal to 0. 
\end{itemize}

\end{rmk}

\begin{ex}
As a specific example, let us use $m=4$, $n=3$. Then, the hypercube with $f(i)=1$ iff $i$ is one of the row vectors of Table \ref{permhcube} is a permutation Costas hypercube; this can be checked by a) verifying that the 2 first columns contain all vectors with integer coordinates between 1 and 3, as do the 2 last columns, and b) by finding all possible $\ds {9\choose 2}=36$ distance vectors and observing they are indeed distinct. Observe also that this hypercube has in total $9=3^2$ nonzero elements out of $81=3^4$; in general, permutation Costas hypercubes have $n^s$ nonzero elements out of $n^{2s}$.
\end{ex}

\begin{table}
\[
\begin{array}{cccc}
		 1 &    1 &    2 &    1\\
     1 &    2 &    2 &    3\\
     1 &    3 &    3 &    1\\
     2 &    1 &    2 &    3\\
     2 &    2 &    1 &    2\\
     2 &    3 &    1 &    3\\
     3 &    1 &    3 &    2\\
     3 &    2 &    3 &    1\\
     3 &    3 &    1 &    2	
\end{array}
\]
\label{permhcube}
\caption{A hypercube with $m=3$, $n=4$ constructed by permuting all pairs of integers}
\end{table}

\begin{rmk} When $m=2s+1,\ s\in\NN$, it is clearly impossible to define a permutation as we did above; the best we can aim for is to find an injective function $g: [n]^s\rightarrow n^{s+1}$, so that $f(i)=1$ iff $(i_{s+1},\ldots,i_{2s+1})= g(i_1,\ldots,i_s),\ (i_1,\ldots,i_s)\in [n]^s$. Obviously, any Costas hypercube in $m-1=2s$ dimensions can be extended to $m$ dimensions, by adding one more coordinate to the position vectors and assigning values randomly to it. This can be done in $n^{s+1}$ ways, as there are $n^s$ vectors and in each one the last (new) coordinate can now take any value in $[n]$: the subvectors formed by the first $2s$ coordinates of the distance vectors remain distinct, so the last coordinate is effectively unused (and therefore purely ornamental). It is of course possible to turn hypercubes in $2s$ dimensions that do not have the Costas property into Costas hypercubes by a judicious choice of the new coordinate's values.
\end{rmk}

\subsection{Strict Costas hypercubes}

Vectors between dots in a permutation Costas hypercube are allowed to have 0 coordinates, in contrast to ordinary 2-dimensional Costas arrays where this does not happen (see Remark \ref{car}). We could potentially restrict the definition of a Costas hypercube to exclude such a possibility, although such a restriction is extremely severe: invoking the Pigeonhole Principle, if a hypercube of side length $n$ has $n+1$ dots, then, fixing any coordinate, at least 2 of those have position vectors with the same value for the chosen coordinate, hence the corresponding distance vector has a 0 there. Therefore, our requirement limits the number of coordinates/dimensions to at most $n$, implying that one of the coordinates (of an element of the hypercube equal to 1) determines unambiguously the rest of them.

\begin{dfn}[Strict Costas hypercubes]
Let $m,n\in\NN$ and let $f: \ZZ^m\rightarrow \{0,1\}$ be such that $f(i)=1$ iff $i\in\{i_1,\ldots,i_n\}$, where $i_j\in[n]^m,\ j\in[m]$, and where all vectors in the family $\{i_j-i_k| 1\leq j<k\leq n\}$ are distinct and have no coordinates equal to 0; then, $f$ will be called a \emph{strict Costas hypercube}. 
\end{dfn}

We proceed to give 2 explicit construction methods for strict Costas hypercubes. 

\begin{thm}[Lifted Costas hypercubes]
Let $m,n\in\NN$ and let $f: \ZZ^m\rightarrow \{0,1\}$. Choose $m-1$ permutations $g_i,\ i\in[m-1]$ (not necessarily distinct), so that $g_1$ is in addition Costas, and set $f(i)=1$ iff $i\in\{(j,g_1(j),\ldots,g_{m-1}(j)):\ j=1,\ldots,n\}$. Then, $f$ is a strict Costas hypercube. 
\end{thm}

\begin{proof}
The result is practically obvious: a typical distance vector is $(j-k,g_1(j)-g_1(k),\ldots,g_{m-1}(j)-g_{m-1}(k))$, where $1\leq j<k\leq n$, and as all $g$s are permutations, no coordinate is equal to 0. Further, the first 2 coordinates of the vectors above are all distinct, as $g_1$ is a Costas permutation. 
\end{proof}

\begin{thm}[Toeplitz-Costas hypercubes]
Let $m,n\in\NN$, $n\leq m$,  and consider the column vector $u=(1,\ldots,n)'$, the vertical shift operator $S(x_1,\ldots,x_n)'=(x_n,x_1,\ldots,x_{n-1})'$, and the array $A=[u\ Su\ \ldots\ S^{m-1}u]$. Let $f: \ZZ^m\rightarrow \{0,1\}$ and set $f(i)=0$ iff $i$ is one of the rows of $A$. Then, $f$ is a strict Costas hypercube and $A$ its \emph{corresponding} Toeplitz array. 
\end{thm}

\begin{proof}
Every column of $A$ contains each integer in $[n]$ exactly once, hence no distance vector has a coordinate equal to 0. Further, given the difference between 2 rows, the 2 rows can be uniquely determined. Indeed, if $r_i$ and $r_j$ denote the rows $i$ and $j$ of $A$, respectively, with $i<j$, then $r_i-r_j$ is a vector with the first $i$ coordinates equal to $i-j$, the following $j-i$ coordinates equal to $n+i-j$, and the remaining coordinates equal to $i-j$; hence, $i$ and $j$ is uniquely determined by inspection. 
\end{proof}

\begin{ex}
Table \ref{toeplitz} shows the corresponding array $A$ of Costas-Toeplitz hypercube with $n=4$ and $m=5$.
\end{ex}

\begin{table}
\[
\begin{array}{ccccc}
		 1 &    4 &    3 &    2 & 1\\
     2 &    1 &    4 &    3 & 2\\
     3 &    2 &    1 &    4 & 3\\
     4 &    3 &    2 &    1 & 4
\end{array}
\]
\label{toeplitz}
\caption{A Toeplitz construction of a Costas hypercube with $m=5$, $n=4$}
\end{table}

The former method has the drawback that it requires that Costas permutations of order $n$ be known, while the latter has the drawback that $m\geq n$, so in general the number of dimensions $m$ needs to be very high; on the other hand, it requires no Costas permutations of order $n$ to be known, so it can be used for values of $n$ such as 32 and 33, where no Costas arrays are known yet. In any case, strict Costas hypercubes have lots of ``empty space'', namely an extremely low density of dots (they only have $n$ dots), and therefore they tend not to be very interesting. 

\section{The main construction method} \label{cn}

It is natural to ask whether Costas arrays (in 2 dimensions) can somehow be manipulated (essentially reshaped) to produce Costas hypercubes, and hopefully permutation Costas hypercubes. We prove below that this is possible, and provide several variants of this construction.   

\subsection{Reshaping}

We formulate and prove below a general result about constructing a Costas hyper-rectangle out of a Costas square. In the special case where the Costas square is a Costas array, and the hyper-rectangle a hypercube, it turns out the hypercube is a permutation Costas hypercube.  

\begin{thm}[Reshaping]\label{resh} Let $m,n\in\NN^*$, $\ds n=\prod_{i=1}^m n_i$, $n_i>1, i\in[m]$, and let $g$ be a Costas permutation of order $n$, but following the convention that $g:[n]-1\rightarrow [n]-1$. Expand $\ds i=\sum_{j=1}^m v_j(i) \prod_{l=j+1}^m n_l$, so that $i$ gets mapped bijectively to $V(i)=(v_1(i),\ldots, v_m(i))$, where $v_j\in[n_j]-1,\ j\in[m]$; similarly, $g(i)$ gets mapped bijectively to $V(g(i))=(v_1(g(i)),\ldots, v_m(g(i)))$. Then, the hyper-rectangle of side length $n_i$ in dimension $i$ and $i+m$, $i\in[m]$, whose dots ($n$ in total) lie at the points $(V(i),V(g(i))):=(v_1(i),\ldots, v_m(i),v_1(g(i)),\ldots, v_m(g(i)))$, $i\in[n]-1$, is actually a permutation Costas hyper-rectangle. 
\end{thm}

\begin{proof} \ 
\begin{enumerate}
	\item Choose 2 values for $i$, say $i_1$ and $i_2$; the corresponding distance vector is:
\begin{multline*}
(V(i_1)-V(i_2),V(g(i_1))-V(g(i_2)))=\\=(v_1(i_1)-v_1(i_2),\ldots, v_m(i_1)-v_m(i_2),v_1(g(i_1))-v_1(g(i_2)),\ldots, v_m(g(i_1))-v_m(g(i_2)))
\end{multline*}
We need to show that all of these vectors are distinct. In other words, we need to show:
\[(V(i_1)-V(i_2),V(g(i_1))-V(g(i_2)))=(V(i_3)-V(i_4),V(g(i_3))-V(g(i_4)))\Rightarrow i_1=i_2,\ i_3=i_4\]
	\item Extend $V^{-1}$ by $\ds V^{-1}(v_1,\ldots,v_m)= \sum_{j=1}^m v_j \prod_{l=j+1}^m n_l$ on the class of vectors where $|v_j|<n_j,\ j\in[m]$. It follows that $V^{-1}(V(i_1)-V(i_2))=i_1-i_2$, as $V(i_1)-V(i_2)$ falls within this class of vectors.
	\item Putting things together:
	\begin{multline*}
	(V(i_1)-V(i_2),V(g(i_1))-V(g(i_2)))=(V(i_3)-V(i_4),V(g(i_3))-V(g(i_4)))\Rightarrow\\
	V^{-1}(V(i_1)-V(i_2),V(g(i_1))-V(g(i_2)))=V^{-1}(V(i_3)-V(i_4),V(g(i_3))-V(g(i_4)))\Leftrightarrow\\ 
	(V^{-1}(V(i_1)-V(i_2)), V^{-1}(V(g(i_1))-V(g(i_2))))=(V^{-1}(V(i_3)-V(i_4)), V^{-1}(V(g(i_3))-V(g(i_4))))\Leftrightarrow \\
	(i_1-i_2, g(i_1)-g(i_2))=(i_3-i_4, g(i_3)-g(i_4))\Leftrightarrow i_1-i_2=i_3-i_4,\ g(i_1)-g(i_2)=g(i_3)-g(i_4) \Rightarrow 
	i_1=i_3,\ i_2=i_4
	\end{multline*}
In the last step we used the fact that $g$ is Costas. 

Further, observe that the left half row vectors so constructed are the expansions of $i\in[n]-1$, while the right half vectors are the expansions of $g(i),\ i\in[n]-1$; the fact that every $i\in[n]-1$ gets expanded exactly once and that $g$ is a permutation guarantees that the hyper-rectangle produced is a permutation one. This completes the proof. 
\end{enumerate}
\end{proof}

The hypercube of even dimension is just a special case:

\begin{cor}[Costas hypercube of even dimension] \label{ch}
Let $m,n\in\NN^*$, and let $g$ be a Costas permutation of order $n^m$, but following the convention that $g:[n]^m-1\rightarrow [n]^m-1$. Expand $\ds i=\sum_{j=1}^m v_j(i) n^{m-j}$, so that $i$ gets mapped bijectively to $V(i)=(v_1(i),\ldots, v_m(i))$, where $v_j\in[n]-1,\ j\in[m]$; similarly, $g(i)$ gets mapped bijectively to $V(g(i))=(v_1(g(i)),\ldots, v_m(g(i)))$. Then, the hypercube of side length $n$ whose dots ($n^m$ in total) lie at the points $(V(i),V(g(i))):=(v_1(i),\ldots, v_m(i),v_1(g(i)),\ldots, v_m(g(i)))$, $i\in[n^m]-1$, is actually a permutation Costas hypercube. 
\end{cor} 

\begin{rmk}
Notice that these hypercubes have $n^m$ dots, which is the square root of the $n^{2m}$ total positions available in the hypercube, just like in Costas arrays, where there are $n$ dots among the $n^2$ available positions. 
\end{rmk}

Here is an attempt to  construct \emph{approximate} Costas hypercubes of odd dimension in the special case where the side length $n$ is a perfect square, using Theorem \ref{resh}, by first constructing a hyper-rectangle as an intermediate step:

\begin{heur}[Costas hypercube of odd dimension]\label{cho}
Let $m,n\in\NN^*$, where $\sqrt{n}\in\NN$, and let $g$ be a Costas permutation of order $n^m\sqrt{n}$, but following the convention that $g:[n]^m\times [\sqrt{n}]-1\rightarrow [n]^m\times [\sqrt{n}]-1$. Expand $\ds i= v_0(i)n^m+\sum_{j=1}^m v_j(i) n^{m-j}$, so that $i$ gets mapped bijectively to $V(i)=(v_0(i),v_1(i),\ldots, v_m(i))$, where $v_j\in[n]-1,\ j\in[m]$ and $v_0\in[\sqrt{n}]-1$; similarly, $g(i)$ gets mapped bijectively to $V(g(i))=(v_0(g(i)), v_1(g(i)),\ldots, v_m(g(i)))$. This process forms a Costas hyper-rectangle in $2m+2$ dimensions, whose side length in $2m$ dimensions is $n$ and in the remaining 2 dimensions $\sqrt{n}$. Now, replace the coordinate pair $(v_0(i), v_0(g(i)))$ in the coordinate vector of each dot by the single coordinate $\sqrt{n}v_0(g(i))+v_0(i)$, $i\in[n^m\sqrt{n}]-1$, which takes values in the range $[n]-1$. Then, the hypercube of side length $n$ whose dots ($n^m\sqrt{n}$ in total) lie at the points $(\sqrt{n}v_0(g(i))+v_0(i),v_1(i),\ldots, v_m(i),v_1(g(i)),\ldots, v_m(g(i)))$, $i\in[n^m\sqrt{n}]-1$, is usually a good approximation of a Costas hypercube, and the removal of a few dots turns it into a Costas hypercube. 
\end{heur}

\begin{rmk} \label{invvf}
The reason why this heuristic often fails to produce a Costas hypercube is that different pairs of coordinates (in difference vectors) can collapse to the same value: for example, assume $n=25\Leftrightarrow \sqrt{n}=5$ and consider the pairs $(-3,0)$ and $(2,-1)$; they get mapped to $-3+5\cdot 0=-3$ and $2-5=-3$. In the context of the proof of Theorem \ref{resh}, this is equivalent to saying that $V(i_1)-V(i_2)$ is \emph{not} always the same as $\text{sign}(i_1-i_2)V(|i_1-i_2|)$. However, simulations show that the damage this does to the Costas property is usually small: tests with Costas arrays of side length $n^3\leq 200$ showed that systematically over 95\% of the difference vectors among the dots are distinct.
\end{rmk}   

The construction methods above can be significantly extended if we use Costas squares instead of Costas arrays as the starting point. The key observation (through the proof of Theorem \ref{resh}) is that the permutation property of the original Costas array is not responsible for the Costas property of the hyper-rectangle produced, but rather for its permutation property alone (and in the case of Heuristic \ref{cho} it does not even achieve that). Therefore, if we are not interested in obtaining a permutation Costas hyper-rectangle as the final product, or if a suitable sized Costas array is not available, we may as well start with a Costas square. 

A special type of Costas squares that proves very helpful in practice is smaller Costas arrays. Consider a Costas array of order $n'\in\NN^*$ and let $n>n'$; then, this Costas array can be turned, by the addition of $n-n'$ blank rows and columns at the sides of the array, into a Costas square of size $n$. Note that such Costas squares are generated by incomplete permutations. We generalize this notion in the following definition:   

\begin{dfn}[Incomplete Costas array] A Costas square with the property that there is at most one dot per row and column will be called an \emph{incomplete Costas array}. A (Costas) hyper-rectangle/hypercube of even dimension with the property that there are no 2 dots whose position vectors have the same left half or right half part will be called an \emph{incomplete (Costas) hyper-rectangle/hypercube}. 
\end{dfn}

\begin{cor}[Constructions out of Costas squares] \ 

\begin{itemize}
	\item A construction according to Theorem \ref{resh} starting with a (incomplete) Costas square results to a (incomplete) Costas hyper-rectangle. 
	\item A construction according to Corollary \ref{ch} starting with a (incomplete) Costas square results to a (incomplete) Costas hypercube.
	\item A construction according to Heuristic \ref{cho} starting with a (incomplete) Costas square results to a hypercube that very nearly has the Costas property and usually can be turned into a Costas hypercube through the removal of a few dots.
\end{itemize}
\end{cor}

\subsection{Construction examples}

We give 3 examples: the first is the construction of a Costas hypercube with $n=5$, $m=4$ out of a Costas array of order 25 (using Corollary \ref{ch}); the second is the construction of a Costas hypercube with $n=9$, $m=3$ out of a Costas array of order 27 (using Heuristic \ref{cho}); and the third is the construction of an incomplete Costas hypercube with $m=5$, $n=4$ (using Heuristic \ref{cho}), starting with a Costas array of order 31 and extending it into an incomplete Costas array of order 32.  

\begin{ex}
Consider the permutation of order 25 appearing on Table \ref{eex} (left). Applying Corollary \ref{ch}, we get a Costas hypercube with $m=4$, $n=5$, whose dot positions appear on Table \ref{eex} (right). Observe that this is a permutation Costas hypercube, as Corollary \ref{ch} states: every vector $(i,j),\ i,j\in[5]-1$ appears in the left 2 columns in exactly 1 row and in the right 2 columns also in exactly one row.  
\end{ex}

\begin{table}
\centering

\begin{tabular}{cc}

\begin{tabular}{|ll|}
\hline
0&10\\\hline
1&7\\\hline
2&6\\\hline
3&9\\\hline
4&17\\\hline
5&23\\\hline
6&21\\\hline
7&2\\\hline
8&20\\\hline
9&11\\\hline
10&15\\\hline
11&3\\\hline
12&12\\\hline
13&22\\\hline
14&1\\\hline
15&16\\\hline
16&18\\\hline
17&19\\\hline
18&5\\\hline
19&0\\\hline
20&13\\\hline
21&24\\\hline
22&14\\\hline
23&8\\\hline
24&4\\\hline
\end{tabular}

&

\begin{tabular}{|llll|}
\hline
0&0&0&2\\\hline
1&0&2&1\\\hline
2&0&1&1\\\hline
3&0&4&1\\\hline
4&0&2&3\\\hline
0&1&3&4\\\hline
1&1&1&4\\\hline
2&1&2&0\\\hline
3&1&0&4\\\hline
4&1&1&2\\\hline
0&2&0&3\\\hline
1&2&3&0\\\hline
2&2&2&2\\\hline
3&2&2&4\\\hline
4&2&1&0\\\hline
0&3&1&3\\\hline
1&3&3&3\\\hline
2&3&4&3\\\hline
3&3&0&1\\\hline
4&3&0&0\\\hline
0&4&3&2\\\hline
1&4&4&4\\\hline
2&4&4&2\\\hline
3&4&3&1\\\hline
4&4&4&0\\\hline
\end{tabular}

\end{tabular}
\caption{\label{eex} The conversion of a Costas array of order 25 into a Costas hypercube with $m=4$, $n=5$: the permutation (left), and the final Costas hypercube (right)}
\end{table}

\begin{ex} 
Consider the Costas permutation of order 27 appearing on Table \ref{oddex} (left). Applying Heuristic \ref{cho}, we first get a Costas hyper-rectangle in 4 dimensions of side lengths 9 and 3, whose dots lie at the points shown in Table \ref{oddex} (center). Subsequently, we combine the middle columns that have values in the range $\{0,1,2\}$ into a single column with values in the range $[9]-1$, as described in Heuristic \ref{cho}: for example, the middle 2 coordinates $(1, 2)$ in row 12 become $1+2\cdot 3=7$. The result appears in Table \ref{oddex} (right). A check of the Costas property shows that this time we get lucky and that the hypercube of $m=3$, $n=9$ we have created is Costas. It is obviously not a permutation Costas hypercube, as this term is meaningless in odd dimensions.   
\end{ex}

\begin{table}
\centering
\begin{tabular}{rrr}

\begin{tabular}{|ll|}
\hline
0&0\\\hline
1&2\\\hline
2&18\\\hline
3&11\\\hline
4&22\\\hline
5&4\\\hline
6&24\\\hline
7&19\\\hline
8&9\\\hline
9&15\\\hline
10&12\\\hline
11&26\\\hline
12&10\\\hline
13&14\\\hline
14&1\\\hline
15&8\\\hline
16&13\\\hline
17&7\\\hline
18&20\\\hline
19&21\\\hline
20&17\\\hline
21&16\\\hline
22&25\\\hline
23&5\\\hline
24&3\\\hline
25&6\\\hline
26&23\\\hline
\end{tabular}

&

\begin{tabular}{|llll|}
\hline
0&0&0&0\\\hline
1&0&0&2\\\hline
2&0&2&0\\\hline
3&0&1&2\\\hline
4&0&2&4\\\hline
5&0&0&4\\\hline
6&0&2&6\\\hline
7&0&2&1\\\hline
8&0&1&0\\\hline
0&1&1&6\\\hline
1&1&1&3\\\hline
2&1&2&8\\\hline
3&1&1&1\\\hline
4&1&1&5\\\hline
5&1&0&1\\\hline
6&1&0&8\\\hline
7&1&1&4\\\hline
8&1&0&7\\\hline
0&2&2&2\\\hline
1&2&2&3\\\hline
2&2&1&8\\\hline
3&2&1&7\\\hline
4&2&2&7\\\hline
5&2&0&5\\\hline
6&2&0&3\\\hline
7&2&0&6\\\hline
8&2&2&5\\\hline
\end{tabular}

&

\begin{tabular}{|lll|}
\hline
0&0&0\\\hline
1&0&2\\\hline
2&6&0\\\hline
3&3&2\\\hline
4&6&4\\\hline
5&0&4\\\hline
6&6&6\\\hline
7&6&1\\\hline
8&3&0\\\hline
0&4&6\\\hline
1&4&3\\\hline
2&7&8\\\hline
3&4&1\\\hline
4&4&5\\\hline
5&1&1\\\hline
6&1&8\\\hline
7&4&4\\\hline
8&1&7\\\hline
0&8&2\\\hline
1&8&3\\\hline
2&5&8\\\hline
3&5&7\\\hline
4&8&7\\\hline
5&2&5\\\hline
6&2&3\\\hline
7&2&6\\\hline
8&8&5\\\hline
\end{tabular}
\end{tabular}
\caption{\label{oddex} The conversion of a Costas array of order 27 into a Costas hypercube with $m=3$, $n=9$: the permutation (left), the intermediate hyper-rectangle (center), and the final Costas hypercube (right)}
\end{table}

\begin{ex}
Consider the Costas permutation of order 31 appearing on Table \ref{oddcs} (left), and consider it as an incomplete Costas permutation of order $32=2^5=4^2\sqrt{4}$. Applying Heuristic \ref{cho}, we first get a Costas hyper-rectangle in 6 dimensions of side lengths 4 and 2, whose dots lie at the points shown in Table \ref{oddcs} (center). Subsequently, we combine the middle columns that have values in the range $\{0,1\}$ into a single column with values in the range $[4]-1$, as described in Heuristic \ref{cho}. The result appears in Table \ref{oddcs} (right). A check of the Costas property shows that this time we get lucky and that the hypercube of $m=5$, $n=4$ we have created is Costas. It is obviously not a permutation Costas hypercube, as this term is meaningless in odd dimensions. It also has just less than $\sqrt{4^5}=32$ dots, one less than that to be exact.   
\end{ex}

\begin{table}

\centering

\begin{tabular}{ccc}

\begin{tabular}{|ll|}
\hline
0&0\\\hline
1&28\\\hline
2&22\\\hline
3&29\\\hline
4&16\\\hline
5&18\\\hline
6&23\\\hline
7&26\\\hline
8&10\\\hline
9&30\\\hline
10&12\\\hline
11&21\\\hline
12&17\\\hline
13&9\\\hline
14&20\\\hline
15&19\\\hline
16&4\\\hline
17&5\\\hline
18&24\\\hline
19&2\\\hline
20&6\\\hline
21&27\\\hline
22&15\\\hline
23&25\\\hline
24&11\\\hline
25&8\\\hline
26&3\\\hline
27&1\\\hline
28&14\\\hline
29&7\\\hline
30&13\\\hline
\end{tabular}

&

\begin{tabular}{|llllll|}
\hline
0&0&0&0&0&0\\\hline
1&0&0&1&3&0\\\hline
2&0&0&1&1&2\\\hline
3&0&0&1&3&1\\\hline
0&1&0&1&0&0\\\hline
1&1&0&1&0&2\\\hline
2&1&0&1&1&3\\\hline
3&1&0&1&2&2\\\hline
0&2&0&0&2&2\\\hline
1&2&0&1&3&2\\\hline
2&2&0&0&3&0\\\hline
3&2&0&1&1&1\\\hline
0&3&0&1&0&1\\\hline
1&3&0&0&2&1\\\hline
2&3&0&1&1&0\\\hline
3&3&0&1&0&3\\\hline
0&0&1&0&1&0\\\hline
1&0&1&0&1&1\\\hline
2&0&1&1&2&0\\\hline
3&0&1&0&0&2\\\hline
0&1&1&0&1&2\\\hline
1&1&1&1&2&3\\\hline
2&1&1&0&3&3\\\hline
3&1&1&1&2&1\\\hline
0&2&1&0&2&3\\\hline
1&2&1&0&2&0\\\hline
2&2&1&0&0&3\\\hline
3&2&1&0&0&1\\\hline
0&3&1&0&3&2\\\hline
1&3&1&0&1&3\\\hline
2&3&1&0&3&1\\\hline
\end{tabular}

&

\begin{tabular}{|lllll|}
\hline
0&0&0&0&0\\\hline
1&0&2&3&0\\\hline
2&0&2&1&2\\\hline
3&0&2&3&1\\\hline
0&1&2&0&0\\\hline
1&1&2&0&2\\\hline
2&1&2&1&3\\\hline
3&1&2&2&2\\\hline
0&2&0&2&2\\\hline
1&2&2&3&2\\\hline
2&2&0&3&0\\\hline
3&2&2&1&1\\\hline
0&3&2&0&1\\\hline
1&3&0&2&1\\\hline
2&3&2&1&0\\\hline
3&3&2&0&3\\\hline
0&0&1&1&0\\\hline
1&0&1&1&1\\\hline
2&0&3&2&0\\\hline
3&0&1&0&2\\\hline
0&1&1&1&2\\\hline
1&1&3&2&3\\\hline
2&1&1&3&3\\\hline
3&1&3&2&1\\\hline
0&2&1&2&3\\\hline
1&2&1&2&0\\\hline
2&2&1&0&3\\\hline
3&2&1&0&1\\\hline
0&3&1&3&2\\\hline
1&3&1&1&3\\\hline
2&3&1&3&1\\\hline
\end{tabular}

\end{tabular}
\caption{\label{oddcs} The conversion of a Costas array of order 31 into a Costas hypercube with $m=5$, $n=4$, treating the Costas array as an incomplete Costas array of order 32: the (incomplete) permutation (left), the intermediate hyper-rectangle (center), and the final Costas hypercube (right)}
\end{table}

\subsection{Applicability of reshaping}

The construction method we described reshapes a Costas array into a high-dimensional Costas hypercube. As the order $n^m$ is bound to be quite big, the only Costas arrays that will normally be available are Golomb and Welch constructions \cite{D}. For practical purposes, we can seek Costas arrays of large orders in databases, such as the database created by Dr. J. K. Beard containing all known Costas arrays up to the order 200 \cite{B}. But more generally, in order to figure out whether a $2m$-dimensional cube of side length $n$ exist, we will need to check whether one of the following equations holds:
\begin{itemize}
	\item $n^m+1=p$ ($W_1$ construction possible)
	\item $n^m+2=p^k$ ($G_2$ construction possible; if $k=1$, $W_2$ construction possible, too)
	\item $n^m+3=p^k$ ($G_3$ construction possible)
\end{itemize}
where in all cases $p$ is a prime. Other variants of the Golomb and the Welch construction do not occur systematically, so we do not investigate them. We now look briefly into each one of the above equations. As a general comment, the solution of these equations falls under the scope of Diophantine Analysis, and it appears that many conjectures can be formulated, but few facts have actually been proved. 

\subsubsection{$n^m+1=p$}

Assume that $m=m_1m_2$, where $m_1$ is odd; then $n^{m_2}+1|n^m+1$, hence it cannot be a prime. It follows that $m=2^k$, for some $k$, whence $p=n^{2^k}+1$. Further, $n$ must necessarily be even, whence $n=2l$ and $p=2^{2^k}l^{2^k}+1$. For $l=1$, we obtain $p=2^{2^k}+1$, the celebrated Fermat primes, for which a lot is known; in particular, it is conjectured that the only such primes correspond to $1\leq k\leq 4$, and lead to hypercubes of side length 2 in $2^{2k}$ dimensions. 

Let us now fix $k=1$, namely seek primes of the form $p=1+n^2$. It is conjectured there are infinitely many such primes \cite{HW}, and it is easy to find several examples: $17=4^2+1,\ 37=6^2+1,\ 101=10^2+1,\ 197=14^2+1,\ 257=16^2+1$ etc. 

Solutions for higher $k$ can also be found. For example, for $1< n\leq 20$ and $k=2$, $n^4+1$ is prime for $n=2,4,6,16,20$. Those solutions can actually build several hypercubes, namely either of side $n^2$ in 4 dimensions, or of side $n$ in 8 dimensions, in all cases with $n^4$ dots. 

\subsubsection{$n^m+2=p^k$}

Here necessarily $n$ must be odd, and again solutions can be found: $3^2+2=11,\ 9^2+2=83,\ 15^2+2=227$; $3^3+2=29,\ 5^3+2=127$; $5^2+2=3^3$ etc. It is not known whether infinitely many solutions exist. This equation appears to be a generalization of Catalan's Equation \cite{Mh}. 

\subsubsection{$n^m+3=p^k$}

Here $n$ must be even, and solutions also exist: $2^4+3=19,\ 2^6+3=67,\ 10^2+3=103,\ 14^2+3=199$ etc. It is not known whether infinitely many solutions exist. This equation also appears to be a generalization of Catalan's Equation \cite{Mh}. 

\subsection{Older generalization attempts}\label{gwgen}

Extensions of the Golomb and the Welch constructions to 3 dimensions were investigated in the past \cite{E}, but the objective then was slightly different: the cubes were so constructed that all the 2-dimensional ``slices'' along some of their directions be Costas arrays, or at least almost Costas arrays in a certain sense. In other words, although the construction was 3-dimensional, the Costas property was still investigated in 2 dimensions. It turns out that many ``reasonable'' generalizations yield cubes with with 2-dimensional (almost) Costas slices. For example, the dots in the cube are placed at the points whose coordinates are the solutions $(i,j,k)$, $i,j,k\in[q-1]$, of one of the following equations:
\begin{enumerate}
	\item $a^{i+x}+b^{j+y}+c^{k+z}=0$, where $a,b,c$ are primitive roots of the field $\FF(q)$, for some $q>2$ power of a prime, and $x,y,z$ fixed in the range $0,\ldots,q-2$;
	\item $a^{i+x}+b^{j+y}\equiv k\mod p$, where $p$ prime, $a,b$ primitive roots of the field $\FF(p)$, and $x,y$ fixed in the range $0,\ldots,q-2$;
	\item $a^{i+j}\equiv k\mod p$, where $p$ prime, and $a$ primitive root of the field $\FF(p)$;
	\item $a^{i+d}\equiv jk\mod p$, where $p$ prime, $a$ primitive root of the field $\FF(p)$, and $d$ fixed in the range $0,\ldots,q-2$.
\end{enumerate}
 
Note that the first one can be expanded in an obvious way into arbitrarily many dimensions. Unfortunately, computer simulations show that none of these constructions yields a Costas cube, so they are not suitable for our purposes. 

\section{An extension of the Welch construction} \label{ewc}

\subsection{The original method and its non-extendability for Costas arrays}

The Welch construction method for Costas arrays \cite{D,G} stipulates that, if $p$ is a prime, $g$ a primitive root \cite{A} of $\FF(p)$ and $c\in[p-1]-1$ a fixed parameter, then the function $\ds f(i)=g^{i-1+c}\mod p$, $i\in[p-1]$ is actually a Costas permutation on $[p-1]$. Contrary to the Golomb construction \cite{D,G}, though, which works in all finite fields, namely with $p^m$ elements where $p$ is a prime and $m\in\NN^*$, the Welch method is not applicable when $m>1$. 

The reason for that is simple: when $m>1$, the elements of the field are no longer represented by integers, but rather by polynomials of degree $m-1$ in a (for all practical purposes) ``dummy'' variable, say $x$, and with coefficients in $[p]-1$, while addition and multiplication are no longer defined modulo an integer, but rather modulo a monic irreducible polynomial $P(x)$ of degree $m$ \cite{A,D}. It follows that the function $f$ of the Welch construction is now $f: [q-1]\rightarrow \FF^*(q)$ where $f(i)=g^{i-1+c}\mod P(x),\ i\in[q-1]$, with $c\in[q-1]-1$ and $g$ a primitive root of $\FF(q)$ (it can be shown that the multiplicative subgroup $\FF^*(q)$ of the field $\FF(q)$ is still cyclic \cite{A,D}); therefore, $f(i)$ is a polynomial, while $i$ is an integer, and the whole construction is (at first sight, at least) meaningless, as we need $f$ to produce integer values! 

Perhaps a more correct way to think of the variable $x$ is that it represents an algebraic element of $\FF(p)$ of order $m$; that is, it is a root of the polynomial $P(x)$ with coefficients in $\FF(p)$ and of degree $m$, while it is the root of no other (non-zero) polynomial with coefficients in $\FF(p)$ and of degree less than $m$ \cite{A}. 

Although this construction fails to yield a Costas array, it still produces Costas hyper-rectangles, as we are about to see.

\subsection{Construction of hyper-rectangles}\label{ewc1}

The field $\FF(p^m)$ can be construed to be a vector space over the field $\FF(p)$ in 2 ways, depending on whether we consider its elements to be polynomials or $m$-tuples:

\begin{dfn}
The field $\FF(p^m)$, where $p$ prime and $m\in\NN^*$, when viewed as a vector space over the field $\FF(p)$, will be denoted by $\FF(p)^m$. Let $V_P(p,m)$ denote the vector space of polynomials of degree $m-1$ over the field $\FF(p)$. Then, $V_P(p,m)$ and $\FF(p)^m$ are isomorphic vector spaces, as they have the same finite dimension $m$ and they are over the same field \cite{A}, under the isomorphism denoted by $\mathcal{F}$, whereby $\mathcal{F}(a_{m-1}x^{m-1}+a_{m-2}x^{m-2}+\ldots+ a_0)= (a_{m-1},a_{m-2},\ldots,a_0)$, $a_i\in\FF(p)$, $i\in[m]-1$.
\end{dfn}

\begin{thm}[Welch hyper-rectangles and hypercubes] \label{welch} Let $p$ be a prime, $m\in\NN^*$, $g$ a primitive root of $\FF(q)$ where $\ds q=p^m$, and $c\in[q-1]-1$. We choose $V_P(p,m)$ as our representation of $\FF(q)$. Then:
\begin{itemize}
	\item The function $f: [q-1]\rightarrow \FF^*(q)$, where $f(i)=\mathcal{F}\left(g^{i-1+c}\mod P(x)\right),\ i\in[q-1]$, $P(x)$ an irreducible polynomial over $\FF(p)$ of degree $m$, is a permutation over $\FF^*(q)$. 
	\item The hyper-rectangle in $m+1$ dimensions with side length $q-1$ in the first dimension and $p$ in the others, whose dots lie at the positions with coordinates $\{(i,f(i))| i\in[q-1]\}$, has the Costas property. 
	\item The hypercube in $2m$ dimensions with side length $p$, whose dots lie at the positions with coordinates $\{(V(i),f(i))| i\in[q-1]\}$, has the Costas property; $V$ is the familiar mapping from Corollary \ref{ch},  with $n=p$ in the present case. 
\end{itemize}
\end{thm}

\begin{proof}
The proof really consists of putting together bits we have already proved. 
\begin{itemize}
	\item $f$ is a permutation over $\FF^*(q)$ because $g$ is taken to be a primitive root of $\FF(q)$. 
	\item That the family of vectors $\{(i,f(i))| i\in[q-1]\}$ has the Costas property follows from a verbatim repetition of the classical argument for $m=1$ \cite{D,G}. Let us see it here in detail: consider the 4 integers $i_1$, $i_2$, $i_3=i_1+k$, and $i_4=i_2+k$ such that $i_j\in[q-1]$, $j\in[4]$, $k\in\NN$, $i_1<i_2$, and let $(i_1-i_2,f(i_1)-f(i_2))=(i_3-i_4,f(i_3)-f(i_4))\Leftrightarrow f(i_1)-f(i_2)=f(i_3)-f(i_4)$. This can be written as:
	\begin{multline*}
	\mathcal{F}\left(g^{i_1-1+c}\mod P(x)\right)-\mathcal{F}\left(g^{i_2-1+c}\mod P(x)\right)=\\=\mathcal{F}\left(g^{i_1+k-1+c}\mod P(x)\right)-\mathcal{F}\left(g^{i_2+k-1+c}\mod P(x)\right)\Rightarrow \text{ ($\mathcal{F}$ is an isomorphism)}\\
	g^{i_1-1+c}-g^{i_2-1+c}\equiv g^{i_1+k-1+c}-g^{i_2+k-1+c}\mod P(x)\Leftrightarrow 
	g^{i_1}-g^{i_2}\equiv g^k(g^{i_1}-g^{i_2})\mod P(x)\Leftrightarrow \\
	(g^k-1)(g^{i_1}-g^{i_2})\equiv 0\mod P(x)\Leftrightarrow g^k-1\equiv 0\mod P(x) \vee g^{i_1}-g^{i_2}\equiv 0\mod P(x)\Leftrightarrow\\
	g^k\equiv 1\mod P(x),\text{ as $0<i_1<i_2<q$, which makes the second condition always false}\Leftrightarrow \\
	k=0,\text{ as $0\leq k\leq q-2$}.
	\end{multline*}
	Hence, $i_1=i_3$, $i_2=i_4$ and the proof is complete.
	\item We need to show that, given that the family $\{(i,f(i))| i\in[q-1]\}$ has the Costas property, the family $\{(V(i),f(i))| i\in[q-1]\}$ has it too; but this is a verbatim repetition of the argument presented in the proof of Theorem \ref{resh}. 
\end{itemize}
\end{proof}

\begin{rmk}\ \label{wl}
\begin{itemize}
	\item The hypercubes constructed above have $q-1$ dots out of $q^2$ possible dot positions, and thus follow approximately the square root rule we saw earlier for the density.
	\item They are ``almost'' permutation hypercubes, except that the zero vectors are missing. If we set $f(0)=0$ and add a dot at the position $(V(0),f(0))=(0,0)$ (obviously $V(0)=0$ as well), then we get a Costas hypercube but we have no guarantee anymore that it has the Costas property; simulations show that sometimes it does. This is the equivalent of the $W_0$ construction of Costas arrays by the addition of a ``corner dot'' to a $W_1$-constructed array \cite{D,G,GT}.
	\item Welch arrays retain the Costas property when their columns get shifted circularly \cite{D,G}; this shift is expressed by the fixed parameter $c$ in the definition of the permutation $f$ shown earlier. As the extension of the Welch method formulated in Theorem \ref{welch} preserves this parameter $c$ in the definition of $f$, we see that the hypercubes and the hyper-rectangles it produces have a similar periodicity: the fact that the family of dots at the locations $\{(i,f(i))| i\in[q-1]\}$ defines a Costas hyper-rectangle implies that the family of dots at the locations $\{(i,f(i\oplus k))| i\in[q-1]\}$, where $i\oplus k:=1+[(i+k-1)\mod (q-1)]$, $k\in[q-1]-1$ fixed, also defines a Costas hyper-rectangle; similarly, the fact that the family of dots at the locations $\{(V(i),f(i))| i\in[q-1]\}$ defines a Costas hypercube implies that the family of dots at the locations $\{(V(i),f(i\oplus k))| i\in[q-1]\}$, $k\in[q-1]-1$ fixed, also defines a Costas hypercube.
	\item Just like we used $V$ to obtain a Costas hypercube out of the originally constructed hyper-rectangle, we can use $V^{-1}$ to obtain a permutation array of order $q-1$, with dots at the positions $\{(i,V^{-1}(f(i)))| i\in[q-1]\}$. The question arises naturally whether this is actually a Costas array; alas, simulations show that it isn't. The reason is that $V^{-1}$ does not generally preserve the Costas property, as we already noted in Remark \ref{invvf}.
	\item The most important aspect of this method is that it builds hypercubes not derived by Costas squares or arrays, at least in an obvious way. At the risk of sounding overly optimistic, if a method that converts Costas hypercubes into Costas arrays were available, it could potentially lead to novel Costas arrays when applied on these hypercubes. 
\end{itemize}
\end{rmk}

\subsection{A further generalization}

The isomorphism $\mathcal{F}$ between the 2 representations of $\FF(p^m)$ as a vector space we suggested in Section \ref{ewc1} is probably the most obvious one, but by no means the only one possible. In terms of basis correspondence, note that $V(p,m)$ is equipped with the natural basis of polynomials $P_i(x)=x^i,\ i\in[m]-1$, $\FF(p)^m$ is equipped with the natural basis of $m$-tuples $e_i=(\delta_{i,1},\ldots,\delta_{i,m-1}),\ i\in[m]-1$, and $\mathcal{F}(P_i)=e_i,\ i\in[m]-1$. Alternatively, we could have used a different isomorphism $\mathcal{F}_B$, such that the polynomial basis $\{P_i|i\in[m]-1\}$ gets mapped to the rows of an invertible matrix $B$ with elements in $\FF(p)$, and consequently: 
\[\mathcal{F}_B(a_{m-1}x^{m-1}+a_{m-2}x^{m-2}+\ldots+ a_0)= (b_{m-1},b_{m-2},\ldots,b_0),\ a_i,b_i\in\FF(p),\ i\in[m]-1\]
where
\[(a_{m-1},a_{m-2},\ldots,a_0)=(b_{m-1},b_{m-2},\ldots,b_0)B\Leftrightarrow (a_{m-1},a_{m-2},\ldots,a_0)B^{-1}=(b_{m-1},b_{m-2},\ldots,b_0)\]
Therefore, $\mathcal{F}_B(\cdot)=\mathcal{F}(\cdot)B^{-1}$. 

Since $B$ is invertible, 
\[(a_{m-1},a_{m-2},\ldots,a_0)=0\Leftrightarrow (b_{m-1},b_{m-2},\ldots,b_0)=0\] 

\begin{thm}[Welch hyper-rectangles and hypercubes under arbitrary bases] \label{welch2} Let $p$ be a prime, $m\in\NN^*$, $g$ a primitive root of $\FF(q)$ where $\ds q=p^m$, $c\in[q-1]-1$, and $B$ an invertible matrix with elements in $\FF(p)$, so that its rows define a basis for the vector space $\FF(p)^m$ over $\FF(p)$. We choose $V_P(p,m)$ as our representation of $\FF(q)$. Then:
\begin{itemize}
	\item The function $f_B: [q-1]\rightarrow \FF^*(q)$, where $f_B(i)=\mathcal{F}(g^{i-1+c}\mod P(x))B^{-1},\ i\in[q-1]$, $P(x)$ an irreducible polynomial over $\FF(p)$ of degree $m$, is a permutation over $\FF^*(q)$. 
	\item The hyper-rectangle in $m+1$ dimensions with side length $q-1$ in the first dimension and $p$ in the others, whose dots lie at the positions with coordinates $\{(i,f_B(i))| i\in[q-1]\}$, has the Costas property. 
	\item The hypercube in $2m$ dimensions with side length $p$, whose dots lie at the positions with coordinates $\{(V(i),f_B(i))| i\in[q-1]\}$, has the Costas property; $V$ is the familiar mapping from Corollary \ref{ch}, with $n=p$ in the present case. 
\end{itemize}
\end{thm}

\begin{proof}\ 

\begin{itemize}
	\item $f_B$ is a permutation over $\FF^*(q)$ because $g$ is taken to be a primitive root of $\FF(q)$ and $\mathcal{F}_B$ is a bijection. 
	\item That the family of vectors $\{(i,f_B(i))| i\in[q-1]\}$ has the Costas property follows from an almost verbatim repetition of the argument presented in Theorem \ref{welch}: consider the 4 integers $i_1$, $i_2$, $i_3=i_1+k$, and $i_4=i_2+k$ such that $i_j\in[q-1]$, $j\in[4]$, $k\in\NN$, $i_1<i_2$, and let $(i_1-i_2,f_B(i_1)-f_B(i_2))=(i_3-i_4,f_B(i_3)-f_B(i_4))\Leftrightarrow f_B(i_1)-f_B(i_2)=f_B(i_3)-f_B(i_4)$. This can be written as:
	\begin{multline*}
	\mathcal{F}_B\left(g^{i_1-1+c}\mod P(x)\right)-\mathcal{F}_B\left(g^{i_2-1+c}\mod P(x)\right)=\\=\mathcal{F}_B\left(g^{i_1+k-1+c}\mod P(x)\right)-\mathcal{F}_B\left(g^{i_2+k-1+c}\mod P(x)\right)\Leftrightarrow \\
	\mathcal{F}\left(g^{i_1-1+c}\mod P(x)\right)B^{-1}-\mathcal{F}\left(g^{i_2-1+c}\mod P(x)\right)B^{-1}=\\=\mathcal{F}\left(g^{i_1+k-1+c}\mod P(x)\right)B^{-1}-\mathcal{F}\left(g^{i_2+k-1+c}\mod P(x)\right)B^{-1}\Leftrightarrow\\
	\mathcal{F}\left(g^{i_1-1+c}\mod P(x)\right)-\mathcal{F}\left(g^{i_2-1+c}\mod P(x)\right)=\\=\mathcal{F}\left(g^{i_1+k-1+c}\mod P(x)\right)-\mathcal{F}\left(g^{i_2+k-1+c}\mod P(x)\right)
	\end{multline*}
	and the rest of the proof proceeds directly as in the proof of Theorem \ref{welch}. 
\end{itemize}
\end{proof}

\begin{rmk}\ 

\begin{itemize}
	\item It is clear that the ability to change the basis while preserving the Costas property increases tremendously the number of possible Welch-constructed Costas hypercubes. 
	\item All these hypercubes are ``almost'' permutation hypercubes: the addition of one dot at the zero position vector turns them into permutation hypercubes.
\end{itemize}
\end{rmk}

\subsection{Normal bases of fields and ``rotational'' constructions}

It is well known that in every finite field $\FF(q)$, $q=p^m$, $p$ prime, $m\in\NN^*$ there exists an element $b$ such that the elements $\ds b^{p^i},\ i=0,\ldots,m-1$ are linearly independent over $\FF(p)$, thus forming a basis of $\FF(p)^m$ \cite{J}. In the context of our generalized Welch construction, these bases (known as \emph{normal} bases) lead to nice ``rotational'' hypercubes. The reason is the equivalence:
\[x=\sum_{i=0}^{m-1}x_i b^{p^i}\Leftrightarrow x^p=\sum_{i=0}^{m-1}x_i b^{p^{i+1}\mod (p^m-1)}=x_{m-1}b+\sum_{i=1}^{m-1}x_{i-1} b^{p^{i}}\]
whereby $x$ has coefficients $(x_0,\ldots,x_{m-1})$ over the normal basis iff $x^p$ has coefficients $(x_{m-1},x_0,\ldots,x_{m-2})$, that is essentially the same expansion only cyclically shifted to the right. 

\subsection{Examples}

\begin{ex}
Let $p=3$ and $m=3$, so that $q=27$, choose $P(x)=x^3+2x+1$ which is irreducible over $\FF(3)$, and choose $c=1$, $g=x$. The Costas hyper-rectangle and the Costas hypercube constructed by Theorem \ref{welch} are shown in Table \ref{twl}, along with the corresponding Welch permutation discussed in Remark \ref{wl}, which fails to have the Costas property. In this particular example, adding a corner dot at $(0,0)$ to the hypercube preserves the Costas property, thus yielding a permutation Costas hypercube.   
\end{ex}

\begin{table}

\centering

\begin{tabular}{ccc}

\begin{tabular}{|llll|}
\hline
1&0&1&0\\\hline
2&1&0&0\\\hline
3&0&1&2\\\hline
4&1&2&0\\\hline
5&2&1&2\\\hline
6&1&1&1\\\hline
7&1&2&2\\\hline
8&2&0&2\\\hline
9&0&1&1\\\hline
10&1&1&0\\\hline
11&1&1&2\\\hline
12&1&0&2\\\hline
13&0&0&2\\\hline
14&0&2&0\\\hline
15&2&0&0\\\hline
16&0&2&1\\\hline
17&2&1&0\\\hline
18&1&2&1\\\hline
19&2&2&2\\\hline
20&2&1&1\\\hline
21&1&0&1\\\hline
22&0&2&2\\\hline
23&2&2&0\\\hline
24&2&2&1\\\hline
25&2&0&1\\\hline
26&0&0&1\\\hline
\end{tabular}

&

\begin{tabular}{|llllll|}
\hline
0&0&1&0&1&0\\\hline
0&0&2&1&0&0\\\hline
0&1&0&0&1&2\\\hline
0&1&1&1&2&0\\\hline
0&1&2&2&1&2\\\hline
0&2&0&1&1&1\\\hline
0&2&1&1&2&2\\\hline
0&2&2&2&0&2\\\hline
1&0&0&0&1&1\\\hline
1&0&1&1&1&0\\\hline
1&0&2&1&1&2\\\hline
1&1&0&1&0&2\\\hline
1&1&1&0&0&2\\\hline
1&1&2&0&2&0\\\hline
1&2&0&2&0&0\\\hline
1&2&1&0&2&1\\\hline
1&2&2&2&1&0\\\hline
2&0&0&1&2&1\\\hline
2&0&1&2&2&2\\\hline
2&0&2&2&1&1\\\hline
2&1&0&1&0&1\\\hline
2&1&1&0&2&2\\\hline
2&1&2&2&2&0\\\hline
2&2&0&2&2&1\\\hline
2&2&1&2&0&1\\\hline
2&2&2&0&0&1\\\hline
\end{tabular}

&

\begin{tabular}{|ll|}
\hline
1&3\\\hline
2&1\\\hline
3&21\\\hline
4&7\\\hline
5&23\\\hline
6&13\\\hline
7&25\\\hline
8&20\\\hline
9&12\\\hline
10&4\\\hline
11&22\\\hline
12&19\\\hline
13&18\\\hline
14&6\\\hline
15&2\\\hline
16&15\\\hline
17&5\\\hline
18&16\\\hline
19&26\\\hline
20&14\\\hline
21&10\\\hline
22&24\\\hline
23&8\\\hline
24&17\\\hline
25&11\\\hline
26&9\\\hline
\end{tabular}

\end{tabular}
\caption{\label{twl} The process of constructing a Welch hypercube in $\FF(27)$: the Welch hyper-rectangle corresponding to $g=x$, $c=0$, $P(x)=x^3+2x+1$ (left), the corresponding Welch hypercube (center), and the corresponding Welch permutation (right)}

\end{table}

\begin{ex}\label{db}
Let $p=5$ and $m=2$, so that $q=25$, $P(x)=x^2+x+2$ which is irreducible over $\FF(5)$, and $c=1$, $g=2x$. Further, choose the array $B$ of Theorem \ref{welch2} to be 
\[
B=
\left[
\begin{array}{cc}
3 & 1\\ 0 & 2
\end{array}
\right]
\Leftrightarrow 
B^{-1}=
\left[
\begin{array}{cc}
3 & 1\\ 0 & 2
\end{array}
\right]
\]
The Costas hyper-rectangle and the Costas hypercube constructed by Theorem \ref{welch2} are shown in Table \ref{twl2}, along with the corresponding Welch permutation discussed in Remark \ref{wl}, which fails to have the Costas property. In this particular example, adding a corner dot at $(0,0)$ to the hypercube preserves the Costas property, thus yielding a permutation Costas hypercube.   
\end{ex}

\begin{table}

\centering

\begin{tabular}{ccc}

\begin{tabular}{|lll|}
\hline
1&4&3\\\hline
2&2&0\\\hline
3&4&1\\\hline
4&1&3\\\hline
5&1&1\\\hline
6&0&4\\\hline
7&2&4\\\hline
8&1&0\\\hline
9&2&3\\\hline
10&3&4\\\hline
11&3&3\\\hline
12&0&2\\\hline
13&1&2\\\hline
14&3&0\\\hline
15&1&4\\\hline
16&4&2\\\hline
17&4&4\\\hline
18&0&1\\\hline
19&3&1\\\hline
20&4&0\\\hline
21&3&2\\\hline
22&2&1\\\hline
23&2&2\\\hline
24&0&3\\\hline
\end{tabular}

&

\begin{tabular}{|llll|}
\hline
0&1&4&3\\\hline
0&2&2&0\\\hline
0&3&4&1\\\hline
0&4&1&3\\\hline
1&0&1&1\\\hline
1&1&0&4\\\hline
1&2&2&4\\\hline
1&3&1&0\\\hline
1&4&2&3\\\hline
2&0&3&4\\\hline
2&1&3&3\\\hline
2&2&0&2\\\hline
2&3&1&2\\\hline
2&4&3&0\\\hline
3&0&1&4\\\hline
3&1&4&2\\\hline
3&2&4&4\\\hline
3&3&0&1\\\hline
3&4&3&1\\\hline
4&0&4&0\\\hline
4&1&3&2\\\hline
4&2&2&1\\\hline
4&3&2&2\\\hline
4&4&0&3\\\hline
\end{tabular}

&

\begin{tabular}{|ll|}
\hline
1&19\\\hline
2&2\\\hline
3&9\\\hline
4&16\\\hline
5&6\\\hline
6&20\\\hline
7&22\\\hline
8&1\\\hline
9&17\\\hline
10&23\\\hline
11&18\\\hline
12&10\\\hline
13&11\\\hline
14&3\\\hline
15&21\\\hline
16&14\\\hline
17&24\\\hline
18&5\\\hline
19&8\\\hline
20&4\\\hline
21&13\\\hline
22&7\\\hline
23&12\\\hline
24&15\\\hline
\end{tabular}

\end{tabular}
\caption{\label{twl2} The process of constructing a Welch hypercube in $\FF(25)$: the Welch hyper-rectangle corresponding to $g=2x$, $c=0$, $P(x)=x^2+x+2$, $B$ chosen as in Example \ref{db} (left), the corresponding Welch hypercube (center), and the corresponding Welch permutation (right)}

\end{table}

\begin{ex}\label{db2}
Let $p=3$ and $m=3$, so that $q=27$, choose $P(x)=x^3+2x+1$ which is irreducible over $\FF(3)$, and choose $c=1$, $g=2x^2$. It follows that $g^3=2x^2+2x+2$, $g^9=2x^2+x+2$, and one can see immediately that the 3 vectors $(g,g^3,g^9)$ are linearly independent. The matrix $B$ corresponding to this choice of a (normal) basis is:

\[
B=\left[
\begin{array}{ccc}
2 & 0 & 0\\
2 & 2 & 2\\
2 & 1 & 2 
\end{array}
\right]
\Leftrightarrow 
B^{-1}=\left[
\begin{array}{ccc}
2 & 0 & 0\\
0 & 1 & 2\\
1 & 1 & 1 
\end{array}
\right]
\]
The Costas hyper-rectangle and the Costas hypercube constructed by Theorem \ref{welch2} are shown in Table \ref{twl3}. In this particular example, adding a corner dot at $(0,0)$ to the hyper-rectangle (either before or after the change of basis) does not preserve the Costas property, but adding it to the hypercube does, thus producing a permutation Costas hypercube.    
\end{ex}

\begin{table}

\centering

\begin{tabular}{ccc}

\begin{tabular}{|llll|}
\hline
1&2&0&0\\\hline
2&1&2&0\\\hline
3&2&2&2\\\hline
4&2&0&2\\\hline
5&2&2&0\\\hline
6&1&0&2\\\hline
7&0&1&0\\\hline
8&0&2&1\\\hline
9&2&1&2\\\hline
10&2&1&1\\\hline
11&0&1&1\\\hline
12&2&2&1\\\hline
13&0&0&2\\\hline
14&1&0&0\\\hline
15&2&1&0\\\hline
16&1&1&1\\\hline
17&1&0&1\\\hline
18&1&1&0\\\hline
19&2&0&1\\\hline
20&0&2&0\\\hline
21&0&1&2\\\hline
22&1&2&1\\\hline
23&1&2&2\\\hline
24&0&2&2\\\hline
25&1&1&2\\\hline
26&0&0&1\\\hline
\end{tabular}

&

\begin{tabular}{|llll|}
\hline
1&0&0&1\\\hline
2&1&2&2\\\hline
3&0&1&0\\\hline
4&2&2&0\\\hline
5&1&2&1\\\hline
6&2&2&1\\\hline
7&2&1&0\\\hline
8&2&0&1\\\hline
9&1&0&0\\\hline
10&0&2&2\\\hline
11&0&2&1\\\hline
12&2&0&2\\\hline
13&2&2&2\\\hline
14&0&0&2\\\hline
15&2&1&1\\\hline
16&0&2&0\\\hline
17&1&1&0\\\hline
18&2&1&2\\\hline
19&1&1&2\\\hline
20&1&2&0\\\hline
21&1&0&2\\\hline
22&2&0&0\\\hline
23&0&1&1\\\hline
24&0&1&2\\\hline
25&1&0&1\\\hline
26&1&1&1\\\hline
\end{tabular}

&

\begin{tabular}{|llllll|}
\hline
0&0&1&0&0&1\\\hline
0&1&0&0&1&0\\\hline
1&0&0&1&0&0\\\hline
0&0&2&1&2&2\\\hline
0&2&0&2&2&1\\\hline
2&0&0&2&1&2\\\hline
0&1&1&2&2&0\\\hline
1&1&0&2&0&2\\\hline
1&0&1&0&2&2\\\hline
0&1&2&1&2&1\\\hline
1&2&0&2&1&1\\\hline
2&0&1&1&1&2\\\hline
0&2&1&2&1&0\\\hline
2&1&0&1&0&2\\\hline
1&0&2&0&2&1\\\hline
0&2&2&2&0&1\\\hline
2&2&0&0&1&2\\\hline
2&0&2&1&2&0\\\hline
1&1&1&2&2&2\\\hline
1&1&2&0&0&2\\\hline
1&2&1&0&2&0\\\hline
2&1&1&2&0&0\\\hline
1&2&2&1&1&0\\\hline
2&2&1&1&0&1\\\hline
2&1&2&0&1&1\\\hline
2&2&2&1&1&1\\\hline
\end{tabular}

\end{tabular}
\caption{\label{twl3} The process of constructing a Welch hypercube in $\FF(27)$: the Welch hyper-rectangle corresponding to $g=2x^2$, $c=1$, $P(x)=x^3+2x+1$ (left), the same hyper-rectangle after changing the basis using $B$ chosen as in Example \ref{db2} (center), and the corresponding Welch hypercube  (right), where we rearranged the order of the rows to demonstrate clearly the rotational structure}

\end{table}


\section{Density of Costas hypercubes}

What is the maximal number of dots we can pack into a hypercube while retaining the Costas property? Up to this point we have avoided on purpose any discussion on this subject, precisely because it is a very complicated one. Note that in the case of Costas arrays this issue is not resolved intrinsically, but rather through the application of an extra condition: we want the array to be a permutation array, hence there can be exactly one dot per row and column. More generally, in a permutation Costas hyper-rectangle there are as many dots as the square root of the total number of elements in the hyper-rectangle. 

This restriction has nothing to do with the Costas property itself: if we lift it (and thus allow any number of dots in a given row or column, even no dots at all, that is), then we can easily find (through simulation) that we can create Costas squares of side $n$ with more than $n$ dots, even for values of $n$ where Costas arrays are not known, like $n=32$ or 33. Similarly, we can create hypercubes with more dots than the permutation restriction would allow, as we show below.

\subsection{Monte Carlo construction of Costas hypercubes}

The determination of the maximum number of dots in a hypercube under the Costas property alone seems to be a difficult combinatorial problem, and we have no solution for it at this time. Even in one dimension, this is the well known problem of the determination of the maximum number of dots in a Golomb ruler of a given length \cite{R}, which still remains unsolved. The best we can do is to try to construct Costas hypercubes through Monte Carlo simulation, using, for example, the algorithm below: 

\begin{alg}\ \label{mca} 

\begin{enumerate}
	\item Set the size $n$ of the hypercube's side and its dimension $m$;
	\item Build a list of all possible $n^m$ coordinate vectors in a random order; 
	\item Place a dot in the hypercube where the first vector in the list indicates;
	\item For each coordinate vector of the list, from the second on to the last, test whether placing a dot in the hypercube where the vector indicates violates the Costas property: if it does, remove it; if it does not, leave it there. 
\end{enumerate}
\end{alg}

We hope that by running this Monte Carlo algorithm enough times we can obtain a Costas hypercube with a density close to the optimal one. Stochastic algorithms on Costas arrays, however, applied on more or less similar problems, are known not to work very well, precisely because Costas arrays are extremely rare \cite{RH}. 

Table \ref{runs} (left) shows the maximum number of dots in Costas hypercubes, for various values of the side length $n$ and the dimension $m$, that we were able to construct by brute force, using the algorithm described above. Some of the hypercubes for which $n=q-2$, $q$ power of a prime, have been constructed by ``sieving'' the Golomb generalization, for comparison purposes: namely, instead of starting with the full list of coordinate vectors in the algorithm we described, we start with the list of vectors satisfying the first Golomb generalization in Section \ref{gwgen}. An important observation here is that using exclusively the latter vectors certainly restricts the potential of the construction method. 

\subsection{Connection with 2 dimensions}

We can apply Algorithm \ref{mca} to construct Costas squares. Such constructions are given in Table \ref{runs} (right), where we see that for small $n$ (approximately $n\leq 50$) the algorithm produces Costas squares more densely packed than Costas arrays, but this seems to be no longer the case for $n>60$. In any case, having produced a Costas square for a given $n$, we can compare it to each Costas array of side length $n$ and find the maximum number of dots that lie in common positions. A systematic high count would show that these arrays are pretty close to Costas arrays, so that, by perturbing the positions of a few dots and possibly removing or adding some others, we could easily obtain a Costas array. Unfortunately, in the tests we ran the count always turned out quite low.

Therefore, in general Costas squares do not help us determine Costas arrays: although they typically have more dots than a Costas array requires, so that it might be expected that by carefully removing some the remaining ones would define a Costas array, in practice they almost invariably contain blank rows and/or columns (without dots, that is), whence it follows that, in addition to dot removal, it is still necessary to move dots around a bit into the blank columns and rows, and the way to do that is not obvious. Another way to state this is that, if we divide the dots of the array into equivalence classes consisting of either the rows or the columns of the array, choosing one representative only from each equivalence class almost invariably leads to incomplete permutations that cannot be completed in an obvious way while retaining the Costas property.

\begin{table}
\[
\begin{array}{|r|r|r|r|} \hline
n & m & \text{Dots} & \text{Golomb}\\ \hline
23 & 3 & 88 & \checkmark \\ \hline
23 & 4 & 100 & \checkmark\\ \hline
47 & 3 & 192 & \checkmark\\ \hline
5 & 3 & 12 & \checkmark\\ \hline
5 & 5 & 72 & \checkmark\\ \hline
17 & 3 & 72 & \checkmark\\ \hline
23 & 3 & 88 & \\ \hline
3 & 5 & 30 & \\ \hline
5 & 3 & 20 & \\ \hline
5 & 5 & 82 & \\ \hline
17 & 3 & 76 & \\ \hline
7 & 4 & 65 & \\ \hline
\end{array}
\hspace{1in}
\begin{array}{|r|r|}\hline
n & \text{Dots}\\ \hline
26 & 32\\ \hline
27 & 33\\ \hline
28 & 34\\ \hline
29 & 35\\ \hline
30 & 36\\ \hline
31 & 36\\ \hline
32 & 37\\ \hline
33 & 38\\ \hline
40 & 43\\ \hline
50 & 52\\ \hline
60 & 58\\ \hline
70 & 65\\ \hline
\end{array}
\]
\caption{\label{runs} The number of dots in some Costas hypercubes constructed by means of Algorithm \ref{mca} (left), as well as in some Costas squares (right). The ticked simulations (left) are produced by running Algorithm \ref{mca} not on the full list of position vectors, but only on those that satisfy the Golomb generalization in Section \ref{gwgen}. }
\end{table}

\section{Summary, conclusion, and future directions}

We defined the multidimensional generalization of Costas arrays in several possible ways, by investigating what the multidimensional generalization of the Costas property should be. We adopted the point of view that the Costas property of a multidimensional binary sequence depends exclusively on its autocorrelation, and that the permutation structure (which we also suitably defined in higher dimensions) is just an extra condition imposed, not directly related to the Costas property itself. 

Hence, Costas arrays became special cases of sequences with the Costas property in 2 dimensions, which we named Costas rectangles and squares; while the former generalized naturally to Costas hyper-rectangles and hypercubes, we generalized the latter in 2 distinct ways:
\begin{itemize}
	\item Strict Costas hypercubes: here, no 2 dots are allowed to have corresponding coordinates with the same value. We proposed 2 construction methods for such hypercubes, but we saw that the requirements of the definition severely limit the number of possible dots, and thus this case tends not to be very interesting. 
	\item Permutation Costas hypercubes: here, the hypercube was viewed as the representation of a permutation between vectors instead of integers. 
\end{itemize}

For the case of permutation Costas hyper-rectangles and hypercubes of even dimension, we proposed a construction method that reshapes an existing Costas array of suitable order into the desired hyper-rectangle or hypercube. In the case of odd dimension, we proposed a related heuristic that does not always work, but even when it doesn't it usually produces a hyper-rectangle that very nearly has the Costas property. We also generalized the constructions by starting with Costas squares instead of arrays, and gave specific examples. 

We subsequently investigated the application of the Welch construction method on finite fields with a nonprime number of elements, and found out that, although the method fails to produce Costas arrays, it produces Costas hyper-rectangles and hypercubes in a natural way, and actually in very large numbers, due to the possibility of changing the basis of the representation while retaining the Costas property. Once more, we supplied specific examples of the construction.

Finally, we investigated experimentally, through Monte Carlo simulations, the difficult question of the restrictions that the Costas property alone imposes on the number of dots in a Costas hypercube, and we generated directly Costas hypercubes and Costas squares of various side lengths in several different dimensions. We observed that, although in small side lengths the simulations produced hypercubes more densely packed with dots than those produced by the construction methods, as the side length increased this ceased to be the case. 

There are still many possible directions for future research in Costas hypercubes. For example,
\begin{enumerate}
	\item What is the maximum number of dots that can be packed into a Costas hypercube or hyper-rectangle, given the number of dimensions and the side lengths? 
	\item Heuristic \ref{cho} seems to be producing Costas hypercubes pretty often, although in many occasions it fails to do so. Can we provide a rigorous and simple sufficient condition for the resulting hypercube to have the Costas property?
	\item Is there a different construction method for Costas hypercubes? In particular, can the Welch and Golomb methods be generalized in a direct way (that is, without the intermediate step of the use of the mapping $V$ defines in Theorem \ref{ch}) in 3 or more dimensions? More generally, can a construction method be found directly based on finite fields?
	\item Are there any engineering applications of Costas hypercubes, perhaps of a similar nature to the applications of Costas arrays?
	\item Is there a method to convert a Costas hypercube into a Costas array? Such a method could potentially lead to the construction of new Costas arrays, because some of the construction methods we have proposed, such as the extended Welch method in Section \ref{ewc}, produce (permutation) Costas hypercubes not linked to any Costas array, at least in an obvious way. 
\end{enumerate}

\section{Acknowledgements}

The author is indebted to Prof. Rod Gow of the Department of Mathematics, University College Dublin, for his suggestion that the construction methods described in Section \ref{cn} could be applied equally well on what we defined as Costas squares, in addition to Costas arrays, as well as for suggesting the change of basis technique in the Welch construction.

\bigskip

\noindent \textbf{Author information:}\\
\\
Konstantinos Drakakis\\
Electronic \& Electrical Engineering\\
University College Dublin\\
Belfield, Dublin 4\\
Ireland\\
\\
Email: Konstantinos.Drakakis@ucd.ie

\end{document}